\documentclass[12pt,a4paper,oneside,reqno,notitlepage]{amsart}
\usepackage{amssymb}
\textwidth
150mm

\theoremstyle{plain}

\numberwithin{equation}{section}

\begin{document}

\baselineskip 8mm
\parindent 9mm

\title[]
{Stochastic delay fractional evolution equations driven by fractional Brownian motion}

\author {Kexue Li}

\address{School of Mathematics and Statistics, Xi'an Jiaotong University, Xi'an 710049, China}

\email{kxli@mail.xjtu.edu.cn}

\thanks{{\it 2000 Mathematics Subjects Classification}: 26A33}

\thanks{This work is partially supported
by the Natural Science Foundation of China under the contact
No. 11201366, China Postdoctoral Science Foundation Funded Project under the contact No. 2012M520080 and Shaanxi Province Postdoctoral Science Foundation Funded Project.}

\begin{abstract}
In this paper, we consider a class of stochastic delay fractional evolution equations driven by fractional Brownian motion in a Hilbert space. Sufficient conditions for the existence and uniqueness of mild solutions are obtained. An application to the stochastic fractional heat equation is presented to illustrate the theory.
\end{abstract}
\maketitle

\section{\textbf{Introduction}}

In recent years, fractional calculus and fractional  equations  have attracted the attention of many researchers due to lots of applications to problems
in physics, chemistry, engineering, biology and finance. One branch of the study is the theory of fractional evolution equations, that is, evolution equations
where the integer-order derivative with respect to time is replaced by a fractional-order derivative. The increasing interest in fractional evolution equations
is motivated by their applications to problems in viscoelasticity, heat conduction in materials with memory, electrodynamics with memory, etc.

Random phenomena exist everywhere in the real world. Systems are often subjected to random perturbations. Stochastic equations have been investigated
by many authors, see, for example, Da Prato and Zabczyk \cite{Da Prato}, Liu \cite{Liu}, van Neerven et al. \cite{Neerven}, Taniguchi \cite{Taniguchi}, Jentzen and R$\ddot{o}$ckner \cite{Arnulf}, Ren et al. \cite{RB}, Ren and Sakthivel \cite{YR}.

There has been some recent interest in studying evolution equations driven by fractional Brownian motion. Duncan et al. \cite{Duncan} studied stochastic differential equations
driven by cylindrical fractional Brownian motion with Hurst parameter  $H\in(\frac{1}{2},1)$. Existence and uniqueness of mild solutions are verified.
Tindel et al. \cite{Tindel} studied linear stochastic evolution equations driven by infinite-dimensional fractional Brownian motion with Hurst parameter in the interval $H\in(0,1)$.
A sufficient and necessary condition for the existence and uniqueness of the
solution is established. Maslowski and Nualart \cite{Maslowski} studied nonlinear stochastic evolution equations in a Hilbert space driven by a cylindrical fractional Brownian motion with Hurst parameter
$H>\frac{1}{2}$ and nuclear covariance operator.  Maslowski and Posp\'{\i}\v{s}il \cite{MJ} proved the existence and egrodicity of a strictly stationary solution for linear stochastic
evolution equations driven by cylindrical fractional Brownian motion.  More recently, Caraballo et al. \cite{Caraballo} investigated the existence, uniqueness and exponential asymptotic bahavior of mild
solutions to a class of stochastic delay evolution equations perturbed by a fractional Brownian motion. Boufoussi and Hajji \cite{Boufoussi} proved the existence and uniqueness of mild solutions of a neutral stochastic differential equations with finite delay, driven by a fractional Brownian motion in a Hilbert space. Ren et al. \cite{YXR} showed the existence and uniqueness of the mild solution for a class of time-dependent stochastic evolution equations with finite delay driven by a standard cylindrical Wiener process and an independent cylindrical fractional Brownian motion.

Some authors have considered  fractional stochastic equations, we refer to EI-Borai \cite{EIBorai}, Ahmed \cite{Ahmed}, Cui and Yan \cite{Cui}, Sakthivel et al. \cite{RPY}, \cite {RPN},  the perturbed terms of these fractional equations are Wiener processes.
To the best of our knowledge, there is no paper which studies the fractional evolution
equations driven by fractional Brownian motion. We study the following
stochastic delay fractional evolution equation
\begin{equation}\label{fbm}
 \ \left\{\begin{aligned} &d[J_{t}^{1-\alpha}(x(t)-x(0))]=Ax(t)dt+f(t,x_{t})dt+h(t)dB^{H}(t), \ t\in [0,b],\\
&x(t)=\phi(t), \ t\in[-r,0],
\end{aligned}\right.
\end{equation}
where $\alpha\in(\frac{1}{2},1]$, $J_{t}^{1-\alpha}$ is the $(1-\alpha)$-order Riemann-Liouville fractional
integral operator, $A$ is the infinitesimal generator of a strongly continuous semigroup $T(t)$, $B^{H}(t)$ is a fractional Brownian motion with Hurst parameter $H\in(\frac{1}{2},1)$, $f,h$ are appropriate functions to be specified later.

The purpose of this paper is to investigate the existence and uniqueness of mild solutions to the equation (\ref{fbm}). The appropriate definition of mild solutions is given and the main results are obtained by using the approximation arguments and an inequality about stochastic integrals with respect to fractional Brownian motion.

The paper is organized as follows. In Section 2, some basic notations and preliminary facts on stochastic integrals for fractional Brownian motion, fractional calculus and some special functions are given. In Section 3, we give the main results of this paper. In Section 4, an example is given to illustrate the obtained theory.

\section{Preliminaries}
In this section we recall some notions and lemmas on Wiener integrals with respect to an infinite dimensional fractional Brownian motion and some basic results about fractional calculus and some special functions.

Let $(\Omega,\mathcal{F},\mathbb{P})$ be a complete probability space. \\
\textbf{Definition 2.1.} An one-dimensional fractional Brownian motion $\{\beta^{H}(t),\ t\geq0\}$ of Hurst index $H\in(0,1)$ is a continuous and centered
Gaussian process with covariance function
\begin{equation*}
R_{H}(t,s)=\mathbb{E}[\beta^{H}(t)\beta^{H}(s)]=\frac{1}{2}(t^{2H}+s^{2H}-|t-s|^{2H}).
\end{equation*}

In the following, we assume $\frac{1}{2}<H<1$.
Consider the square integrable kernel
\begin{equation*}
K_{H}(t,s)=c_{H}s^{\frac{1}{2}-H}\int_{s}^{t}(u-s)u^{H-\frac{1}{2}}du,
\end{equation*}
where $c_{H}=\left[\frac{H(2H-1)}{\beta(2-2H,H-\frac{1}{2}))}\right]^{1/2}$, $t>s$.\\
Then
\begin{equation*}
\frac{\partial K_{H}}{\partial t}(t,s)=c_{H}\left(\frac{t}{s}\right)^{H-\frac{1}{2}}(t-s)^{H-\frac{3}{2}}.
\end{equation*}

Consider a fBM $\{\beta^{H}(t), \ t\in[0,b]\}$. We denote by $\zeta$ the set of step functions on $[0,b]$. Let $\mathcal{H}$ be the Hilbert space defined as
the closure of $\zeta$ with respect to the scalar product
\begin{equation*}
\langle \textbf{1}_{[0,t]}, \textbf{1}_{[0,s]}\rangle_{\mathcal{H}}=R_{H}(t,s).
\end{equation*}
The mapping $\textbf{1}_{[0,t]}\rightarrow \{\beta^{H}(t)\}$ can be extended to an isometry between $\mathcal{H}$ and the first
Wiener chaos of the fBM $\overline{\mbox{span}}^{L^{2}(\Omega)}\{\beta^{H}(t), \ t\in[0,b]\}$.

Consider the linear operator $K_{H}^{*}$ from  $\zeta$ to $L^{2}([0,b])$ defined by
\begin{equation*}
(K_{H}^{*}\varphi)(s)=\int_{s}^{b}\varphi(t)\frac{\partial K_{H}}{\partial t}(t,s)dt.
\end{equation*}

The operator $K_{H}^{*}$  is an isometry between $\zeta$ and $L^{2}([0,b])$ that can be extended to the Hilbert space $\mathcal{H}$.

Consider the process $W=W(t),t\in[0,b]$ defined by
\begin{equation*}
W(t)=\beta^{H}((K_{H}^{\ast})^{-1}\textbf{1}_{[0,t]}).
\end{equation*}
Then $W$ is a Wiener process and $\beta^{H}$ has the integral representation
\begin{equation*}
\beta^{H}(t)=\int_{0}^{t}K_{H}(t,s)dW(s).
\end{equation*}
Let $|\mathcal{H}|$ be the Banach space of measurable functions $\varphi$ on $[0,b]$ such that
\begin{equation*}
\|\varphi\|_{|\mathcal{H}|}^{2}=\alpha_{H}\int_{0}^{b}\int_{0}^{b}\varphi(r)\varphi(u)|r-u|^{2H-2}dudr<\infty,
\end{equation*}
where $\alpha_{H}=H(2H-1)$.

We have the embeddings (see \cite{Nualart})
\begin{equation*}
L^{1}([0,t])\subset L^{\frac{1}{H}} \subset |\mathcal{H}|\subset\mathcal{H}.
\end{equation*}
Throughout this paper, $X$ and $Y$ denote two real separable Hilbert spaces. By $\mathcal{L}(Y,X)$, we denote the space of all bounded linear operators from $Y$ to $X$. For convenience, we use the notation $\|\cdot\|$ to denote
the norms in $X, Y$ and $\mathcal{L}(Y,Y)$ when no confusion possibly arises. Suppose that there exists a complete orthonormal system $\{e_{n}\}_{n=1}^{\infty}$ in $Y$, let $Q\in \mathcal{L}(Y,Y)$ be an operator
with finite trace $trQ=\sum_{n=1}^{\infty}\lambda_{n}<\infty \ (\lambda_{n}\geq 0)$ such that $Qe_{n}=\lambda_{n}e_{n}$. The infinite dimensional fBM on $Y$ can be defined by using
covariance operator $Q$ as
\begin{equation*}
B^{H}(t)=B_{Q}^{H}(t)=\sum_{n=1}^{\infty}\sqrt{\lambda_{n}}e_{n}\beta_{n}^{H}(t),
\end{equation*}
where $\beta_{n}^{H}(t)$ are one dimensional standard fractional Brownian motions mutually independent on $(\Omega,\mathcal{F},\mathbb{P})$.
Consider the space $\mathcal{L}_{2}^{0}:=\mathcal{L}_{2}^{0}(Y,X)$ of all $Q$-Hilbert-Schmidt operators $\varphi: Y\rightarrow X$. We recall that $\varphi\in \mathcal{L}(Y,X)$ is called
a $Q$-Hilbert-Schmidt operator, if
\begin{equation*}
\|\varphi\|_{\mathcal{L}_{2}^{0}}^{2}:=\sum_{n=1}^{\infty}\|\sqrt{\lambda_{n}}\varphi e_{n}\|^{2}<\infty,
\end{equation*}
and that the space $\mathcal{L}_{2}^{0}$ equipped with the inner product $\langle\varphi,\psi\rangle_{\mathcal{L}_{2}^{0}}=\sum_{n=1}^{\infty}\langle\varphi e_{n}, \psi e_{n}\rangle$ is a separable Hilbert space.

Let $(\phi(s))_{s\in[0,b]}$ be a deterministic function with values in $\mathcal{L}_{2}^{0}(Y,X)$. The stochastic integral of $\phi$ with respect to $B^{H}$ is defined by
\begin{equation}\label{definition}
\int_{0}^{t}\phi(s)dB^{H}(s)=\sum_{n=1}^{\infty}\int_{0}^{t}\sqrt{\lambda_{n}}\phi(s)e_{n}d\beta_{n}^{H}(s)=\sum_{n=1}^{\infty}\int_{0}^{t}\sqrt{\lambda_{n}}(K_{H}^{*}(\phi e_{n}))(s)d\beta_{n}(s).
\end{equation}
\textbf{Lemma 2.1} (\cite{Boufoussi}). If $\varphi: [0,b]\rightarrow \mathcal{L}_{2}^{0}(Y,X)$ satisfies $\int_{0}^{b}\|\varphi(s)\|_{\mathcal{L}_{2}^{0}}^{2}<\infty$ then the above sum in
(\ref{definition}) is well defined as an $X$-valued random variable and we have
\begin{equation*}
\mathbb{E}\|\int_{0}^{t}\varphi(s)dB^{H}(s)\|^{2}\leq 2Ht^{2H-1}\int_{0}^{t}\|\varphi(s)\|_{\mathcal{L}_{2}^{0}}^{2}ds.
\end{equation*}

Now we recall some notations and preliminary results about fractional calculus and some special functions.

For $\beta>0$, let
\begin{eqnarray}
g_\beta(t)=\left\{\begin{aligned}
&\frac{t^{\beta-1}}{\Gamma(\beta)},\;&
t>0,\\
&0,\; &t\leq 0,
\end{aligned}\right.
\end{eqnarray}
where $\Gamma(\cdot)$ is the Gamma function. Let
$g_{0}(t)=\delta(t)$, the delta distribution.\\
\textbf{Definition 2.1.} The Riemann-Liouville fractional integral
of order $\alpha>0$ of $f:[0,b]\rightarrow X$ is defined by
\begin{equation}\label{0}
J_{t}^{\alpha}f(t)=\frac{1}{\Gamma(\alpha)}\int_{0}^{t}(t-s)^{\alpha-1}f(s)ds.
\end{equation}
\textbf{Definition 2.2.} The Riemann-Liouville fractional derivative
of order $\alpha\in (0,1]$ of $f:[0,b]\rightarrow X$ is defined by
\begin{equation}\label{ft}
D_{t}^{\alpha}f(t)=\frac{d}{dt}J_t^{1-\alpha}f(t).
\end{equation}
\textbf{Definition 2.3.} The Caputo fractional derivative
of order $\alpha\in (0,1]$ of $f:[0,b]\rightarrow X$ is defined by
\begin{equation}\label{dt}
^{C}D_{t}^{\alpha}f(t)=D_{t}^{\alpha}(f(t)-f(0)).
\end{equation}
The Laplace transform of the Caputo fractional
derivative is given by
\begin{equation}\label{laplace}
L\{^{C}D_{t}^{\alpha}u(t)\}=\lambda^{\alpha}\hat{u}(\lambda)-\lambda^{\alpha-1}u(0),
\end{equation}
where $\hat{u}(\lambda)$ is the Laplace transform of $u$ defined by
\begin{equation}\label{cute}
\hat{u}(\lambda)=\int_{0}^{\infty}e^{-\lambda t}u(t)dt,\
\mathfrak{R}\lambda>\omega,
\end{equation}
where $\mathfrak{R}\lambda$ stands for the real part of the complex number
$\lambda$.\\
\textbf{Definition 2.4}  The Mittag-Leffler function
 is defined by
\begin{equation}\label{ML}
E_{\alpha,\beta}(z)=\sum_{n=0}^{\infty}\frac{z^{n}}{\Gamma(\alpha
n+\beta)},\ \alpha,\beta>0,\ z\in \mathbb{C}.
\end{equation}
When $\beta=1$, set $E_{\alpha}(z)=E_{\alpha,1}(z)$. \\
\textbf{Definition 2.5} (\cite{Igor}). The Mainardi's function is defined by
\begin{equation}\label{rep2}
M_{\alpha}(z)=\sum_{n=0}^{\infty}\frac{(-z)^{n}}{n!\Gamma(-\alpha n+1-\alpha)},\ 0<\alpha<1,\ z\in \mathbb{C}.
\end{equation}
The Laplace transform of the Mainardi's function $M_{\alpha}(r)$ is (\cite{Main}):
\begin{align}\label{lap}
\int_{0}^{\infty}e^{-r\lambda}M_{\alpha}(r)dr=E_{\alpha}(-\lambda).
\end{align}
By (\ref{ML}) and (\ref{lap}), it is clear that
\begin{align}\label{one}
\int_{0}^{\infty}M_{\alpha}(r)dr=1,\ 0<\alpha<1.
\end{align}
On the other hand,
$M_{\alpha}(z)$ satisfies
the following equality (\cite{Main}, P.249)
\begin{equation}\label{L}
\int_{0}^{\infty}\frac{\alpha}{r^{\alpha+1}}M_{\alpha}(1/r^{\alpha})e^{-\lambda r}dr=e^{-\lambda^{\alpha}}.
\end{equation}
and the equality (\cite{Main}, (F.33))
\begin{align}\label{rm}
\int_{0}^{\infty}r^{\delta}M_{\alpha}(r)dr=\frac{\Gamma(\delta+1)}{\Gamma(\alpha\delta+1)}, \ \delta>-1,\ 0< \alpha<1.
\end{align}

\section{Main results}
Denote by $C([a,b]; L^{2}(\Omega;X))=C([a,b]; L^{2}(\Omega,\mathcal{F},\mathbb{P};X))$ the Banach space of all continuous functions from $[a,b]$ into $L^{2}(\Omega;X)$
equipped with the sup norm. Let $r>0$ be a real number. If $x\in C([-r,b];L^{2}(\Omega;X))$ for each $t\in [0,b]$, we denote by $x_{t}\in C([-r,b];L^{2}(\Omega;X))$ the function defined by $x_{t}(s)=x(t+s)$ for $s\in [-r,0]$.

Note that if $g=0$ and the delay term in (\ref{fbm}) is omitted, thanks to (\ref{dt}) and (\ref{ft}), we can consider the following deterministic case:
\begin{equation}\label{rb}
\left \{\begin{aligned}
&^{C}D_{t}^{\alpha}u(t)=Au(t)+f(t),\ t\in [0,b],\\
& u(0)=\eta,
\end{aligned}\right.
\end{equation}
We refer to \cite{kexue} about the definition of mild solutions to (\ref{rb}), for convenience, the details is presented in the following.
Since the coefficient operator $A$ is the infinitesimal generator of a $C_{0}$-semigroup $T(t)$, there exists constants $N\geq 1,\ \omega\geq 0$, such that $\|T(t)\|\leq Ne^{\omega t}$,  \ $t\geq 0$.
From the subordination principle (see Theorem 3.1 in \cite{Bazhlekova}), it follows that
$A$ generates an exponentially bounded solution operator $T_{\alpha}(t)(0<\alpha<1)$ (see \cite{Bazhlekova}, Definition 2.3), satisfying $T_{\alpha}(0)=I$ (the identity operator on $X$),
\begin{align}\label{exp}
\|T_{\alpha}(t)\|\leq Ne^{\omega^{1/\alpha}t}, \ t\geq0,
\end{align}
and
\begin{align}\label{rep1}
T_{\alpha}(t)&=\int_{0}^{\infty}t^{-\alpha}M_{\alpha}(st^{-\alpha})T(s)ds\nonumber\\
&=\int_{0}^{\infty}M_{\alpha}(r)T(t^{\alpha}r)dr,\ t>0,
\end{align}
where $M_{\alpha}(r)$ is the Mainardi's function.

Since $\|T_{\alpha}(t)\|\leq N_{1}e^{\omega^{1/\alpha}t}$, from the formulas (2.5) and (2.6) in \cite{Bazhlekova}, it follows that
\begin{equation}\label{rep3}
\{\lambda^{\alpha}:\lambda>\omega^{1/\alpha}\}\subset \rho(A),
\end{equation}
and
\begin{equation}\label{rep4}
\lambda^{\alpha-1}R(\lambda^{\alpha},A)\eta=\int_{0}^{\infty}e^{-\lambda t}T_{\alpha}(t)\eta dt,\ \lambda>\omega^{1/\alpha},\ \eta\in X.
\end{equation}
By $\widehat{f}(\lambda)=\int_{0}^{\infty}e^{-\lambda t}f(t)dt$, resp. $\widehat{u}(\lambda)=\int_{0}^{\infty}e^{-\lambda t}u(t)dt$, we denote the Laplace transforms of the functions $f\in L^{1}([0,b];X)$ and $u\in C([0,b];X)$, respectively.
Taking the Laplace transform to both sides of $^{C}D_{t}^{\alpha}u(t)=Au(t)+f(t)$ and using the initial condition $u(0)=\eta$, we obtain
\begin{equation}\label{rep5}
\widehat{u}(\lambda)=\lambda^{\alpha-1}R(\lambda^{\alpha},A)\eta+R(\lambda^{\alpha},A)\widehat{f}(\lambda).
\end{equation}
Since $A$ is the infinitesimal generator of $C_{0}$-semigroup $T(t)$,
\begin{equation}\label{rep6}
R(\lambda,A)\eta=\int_{0}^{\infty}e^{-\lambda t}T(t)\eta dt,\ \lambda>\omega,\ \eta\in X.
\end{equation}
Hence, by (\ref{rep6}) and (\ref{L}),
\begin{align}\label{transfer}
&R(\lambda^{\alpha},A)\widehat{f}(\lambda)\\
&=\int_{0}^{\infty}e^{-\lambda^{\alpha}t}T(t)\widehat{f}(\lambda)dt\nonumber\\
&=\int_{0}^{\infty}\int_{0}^{\infty}\alpha t^{\alpha-1}e^{-(\lambda t)^{\alpha}}T(t^{\alpha})f(s)e^{-s\lambda}dsdt\nonumber\\
&=\int_{0}^{\infty}\int_{0}^{\infty}\int_{0}^{\infty}\alpha t^{\alpha-1}\frac{\alpha}{r^{\alpha+1}}M_{\alpha}(1/r^{\alpha})e^{-\lambda tr}T(t^{\alpha})f(s)e^{-s\lambda}drdsdt\nonumber\\
&=\int_{0}^{\infty}\int_{0}^{\infty}\int_{0}^{\infty}\frac{\alpha^{2}}{r^{2\alpha+1}}t^{\alpha-1}M_{\alpha}(1/r^{\alpha})e^{-\lambda t}f(s)e^{-s\lambda}drdsdt\nonumber\\
&=\int_{0}^{\infty}\int_{0}^{\infty}\int_{0}^{\infty}\alpha rt^{\alpha-1}M_{\alpha}(r)T(t^{\alpha}r)f(s)e^{-\lambda(t+s)}drdsdt\nonumber\\
&=\int_{0}^{\infty}e^{-\lambda t}\left(\int_{0}^{t}\int_{0}^{\infty}\alpha rM_{\alpha}(r)(t-s)^{\alpha-1}T((t-s)^{\alpha}r)f(s)drds\right)dt.
\end{align}
By (\ref{rep4}), (\ref{rep5}), (\ref{transfer}), we obtain
\begin{align}\label{trans}
\widehat{u}(\lambda)&=\int_{0}^{\infty}e^{-\lambda t}T_{\alpha}(t)\eta dt\nonumber\\
&\quad+\int_{0}^{\infty}e^{-\lambda t}\left(\int_{0}^{t}\int_{0}^{\infty}\alpha rM_{\alpha}(r)(t-s)^{\alpha-1}T((t-s)^{\alpha}r)f(s)drds\right)dt.
\end{align}
From the uniqueness theorem of the Laplace transform, it follows that
\begin{align}\label{tran}
u(t)=T_{\alpha}(t)\eta+\int_{0}^{t}\int_{0}^{\infty}\alpha rM_{\alpha}(r)(t-s)^{\alpha-1}T((t-s)^{\alpha}r)f(s)drds.
\end{align}
Set
\begin{align}\label{t}
S_{\alpha}(t)\eta=\int_{0}^{\infty}\alpha rM_{\alpha}(r)T(t^{\alpha}r)\eta dr, \ t\geq 0, \ \eta\in X.
\end{align}
By (\ref{tran}), (\ref{t}),
\begin{align}\label{tr}
u(t)=T_{\alpha}(t)\eta+\int_{0}^{t}(t-s)^{\alpha-1}S_{\alpha}(t-s)f(s)ds,
\end{align}
where
\begin{align}\label{rep1}
T_{\alpha}(t)=\int_{0}^{\infty}M_{\alpha}(r)T(t^{\alpha}r)dr,\quad t\geq0,
\end{align}
\begin{align}\label{tra}
S_{\alpha}(t)=\int_{0}^{\infty}\alpha rM_{\alpha}(r)T(t^{\alpha}r)dr, \quad t\geq 0,
\end{align}
where $M_{\alpha}(r)$ is the Mainardi's function, $T(t)$ is the semigroup generated by $A$.\\

Motivated by (\ref{tr}), we give the following definition of mild solutions of Eq. (\ref{fbm}).\\
\textbf{Definition 3.1.} A $X$-valued process $x(t)$ is called a mild solution of Eq. (\ref{fbm}) if\\
(i) $x(\cdot)\in C([-r,b];L^{2}(\Omega;X))$,\\
(ii) $x(t)=\phi(t)$,\ $t\in [-r,0]$,\\
(iii) For all $t\in[0,b]$,
\begin{align}\label{solution}
x(t)&=T_{\alpha}(t)\phi(0)+\int_{0}^{t}(t-s)^{\alpha-1}S_{\alpha}(t-s)f(s,x_{s})ds\nonumber\\
&\quad+\int_{0}^{t}(t-s)^{\alpha-1}S_{\alpha}(t-s)h(s)dB^{H}(s) \ \ \mathbb{P}-a.s.
\end{align}

In order to study the existence and uniqueness of mild solutions of Eq. (\ref{fbm}), we make the following assumptions:\\
$(a_{1})$ $A: D(A)\subset X\rightarrow X$ is the infinitesimal generator of a strongly continuous semigroup of bounded linear operators $T(t)$ on $X$, there exists constants $M\geq 1$, $\omega\geq 0$ such that $\|T(t)\|\leq Me^{\omega t}.$\\
The function $f:[0,b]\times C([-r,0];X)\rightarrow X$  satisfies the following conditions:\\
$(f_{1})$ The mapping $t\in(0,b)\rightarrow f(t,\xi)\in X$ is Lebesgue measurable, for a.e. $t$ and for all $\xi\in C([-r,0];X)$.\\
$(f_{2})$ There exists a positive constant $C_{f}>0$ such that for any  $x,\ y\in C([-r,0];X)$ and $t\in [0,b]$,
\begin{align*}
\int_{0}^{t}\|f(s,x_{s})-f(s,y_{s})\|^{2}ds\leq C_{f}\int_{-r}^{t}\|x(s)-y(s)\|^{2}ds.
\end{align*}
$(f_{3})$
\begin{align*}
\int_{0}^{b}\|f(s,0)\|^{2}ds<\infty.
\end{align*}
$(h_{1})$ There exists a constant $p>\frac{1}{2\alpha-1}$ such that the function $h: [0,\infty)\rightarrow L_{2}^{0}(Y,X)$ satisfies
\begin{align*}
\int_{0}^{b}\|h(s)\|_{L_{2}^{0}}^{2p}ds<\infty.
\end{align*}
\textbf{Remark 3.2.}  If $\|T(t)\|\leq Me^{\omega t}$ for $M\geq 1$, $\omega\geq 0$, then
$\|S_{\alpha}(t)\|\leq ME_{\alpha,\alpha}(\omega t^{\alpha})$. \\
\textbf{Proof.} By (\ref{tra}), (\ref{rm}) and the property of the Gamma function,  we have
\begin{align*}
\|S_{\alpha}(t)\|&\leq \alpha M\int_{0}^{\infty}rM_{\alpha}(r)e^{\omega t^{\alpha}r}dr\\
&\leq \alpha M\int_{0}^{\infty}rM_{\alpha}(r)\sum_{n=0}^{\infty}\frac{(\omega t^{\alpha}r)^{n}}{n!}dr\\
&\leq \alpha M\sum_{n=0}^{\infty}\frac{(\omega t^{\alpha})^{n}}{n!}\int_{0}^{\infty}r^{n+1}M_{\alpha}(r)dr\\
&\leq \alpha M\sum_{n=0}^{\infty}\frac{(\omega t^{\alpha})^{n}}{n!}\cdot \frac{\Gamma(n+2)}{\Gamma(\alpha(n+1)+1)}\\
&\leq \alpha M\sum_{n=0}^{\infty}\frac{(\omega t^{\alpha})^{n}}{n!}\cdot \frac{(n+1)\Gamma(n+1)}{\alpha(n+1)\Gamma(\alpha(n+1))}\\
&\leq M E_{\alpha,\alpha}(\omega t^{\alpha}).
\end{align*}
\textbf{Theorem 3.3.} Suppose that $(a_{1})$, $(f_{1})$, $(f_{2})$, $(f_{3})$ and $(h_{1})$ are satisfied. Then for every $\phi \in C([-r,0];L^{2}(\Omega;X))$, Eq.(\ref{fbm}) has a unique mild solution on $[-r,b]$.\\
\textbf{Proof.} First, we prove the uniqueness of mild solutions. Suppose that $x,y\in C([-r,b];L^{2}(\Omega;X))$ are two mild solutions of (\ref{fbm}).
\begin{align}\label{gron}
\mathbb{E}\|x(t)-y(t)\|^{2}&=\mathbb{E}\|\int_{0}^{t}(t-s)^{\alpha-1}S_{\alpha}(t-s)(f(s,x_{s})-f(s,y_{s}))ds\|^{2}\nonumber\\
&\leq M^{2}(\mathbb{E}_{\alpha,\alpha}(\omega t^{\alpha}))^{2} \mathbb{E}\left(\int_{0}^{t}(t-s)^{\alpha-1}\|f(s,x_{s})-f(s,y_{s})\|ds\right)^{2}\nonumber\\
&\leq M^{2}(E_{\alpha,\alpha}(\omega t^{\alpha}))^{2} \left(\int_{0}^{t}(t-s)^{2\alpha-2}ds\right)\mathbb{E}\int_{0}^{t}\|f(s,x_{s})-f(s,y_{s})\|^{2}ds\nonumber\\
&\leq \frac{M^{2}(E_{\alpha,\alpha}(\omega t^{\alpha}))^{2}t^{2\alpha-1}C_{f}}{(2\alpha-1)}\int_{0}^{t}\mathbb{E}\|x(s)-y(s)\|^{2}ds\nonumber\\
&\leq \frac{M^{2}(E_{\alpha,\alpha}(\omega t^{\alpha}))^{2}t^{2\alpha-1}C_{f}}{(2\alpha-1)}\int_{0}^{t}\sup_{0\leq\tau\leq s}\mathbb{E}\|x(\tau)-y(\tau)\|^{2}ds.
\end{align}
Taking the supremum in (\ref{gron}), we have
\begin{align}\label{uni}
\sup_{0\leq \mu \leq t}\mathbb{E}\|x(\mu)-y(\mu)\|^{2}ds\leq \frac{M^{2}(E_{\alpha,\alpha}(\omega b^{\alpha}))^{2}b^{2\alpha-1}C_{f}}{(2\alpha-1)} \int_{0}^{t}\sup_{0\leq\tau\leq s}\mathbb{E}\|x(\tau)-y(\tau)\|^{2}ds.
\end{align}
By the Gronwall inequality, we have
\begin{align}
\mathbb{E}\|x(t)-y(t)\|^{2}=0,\ t\in[0,b].
\end{align}
This together with $x(t)=y(t)$ on $[-r,0]$ yield that the mild solution is unique.

Next we prove the existence of mild solutions of Eq.(\ref{fbm}).\\
Step 1: We show that the stochastic integral term $\int_{0}^{t}(t-s)^{\alpha-1}S_{\alpha}(t-s)h(s)dB^{H}(s)$ possesses the required regularity. For sufficiently small $\delta>0$,
\begin{align}\label{delta}
&\mathbb{E}\|\int_{0}^{t+\delta}(t+\delta-s)^{\alpha-1}S_{\alpha}(t+\delta-s)h(s)dB^{H}(s)-\int_{0}^{t}(t-s)^{\alpha-1}S_{\alpha}(t-s)h(s)dB^{H}(s)\|^{2}\nonumber\\
&\leq 2\mathbb{E}\|\int_{0}^{t}((t+\delta-s)^{\alpha-1}S_{\alpha}(t+\delta-s)-(t-s)^{\alpha-1}S_{\alpha}(t-s))h(s)dB^{H}(s)\|^{2}\nonumber\\
&\quad+2\mathbb{E}\|\int_{t}^{t+\delta}(t+\delta-s)^{\alpha-1}S_{\alpha}(t+\delta-s)h(s)dB^{H}(s)\|^{2}\nonumber\\
&=I_{1}+I_{2}.
\end{align}
For $I_{1}$, by Lemma 2.1, we have
\begin{align}\label{I1}
I_{1}&=2\mathbb{E}\|\int_{0}^{t}((t+\delta-s)^{\alpha-1}-(t-s)^{\alpha-1})S_{\alpha}(t+\delta-s)h(s)dB^{H}(s)\nonumber\\
&\quad+\int_{0}^{t}(t-s)^{\alpha-1}(S_{\alpha}(t+\delta-s)-S_{\alpha}(t-s))g(s)dB^{H}(s)\|^{2}\nonumber\\
&\leq 4\mathbb{E}\|\int_{0}^{t}((t+\delta-s)^{\alpha-1}-(t-s)^{\alpha-1})S_{\alpha}(t+\delta-s)h(s)dB^{H}(s)\|^{2}\nonumber\\
&\quad+4\mathbb{E}\|\int_{0}^{t}(t-s)^{\alpha-1}(S_{\alpha}(t+\delta-s)-S_{\alpha}(t-s))h(s)dB^{H}(s)\|^{2}\nonumber\\
&\leq8Ht^{2H-1}M^{2}(E_{\alpha,\alpha}(\omega (t+\delta)^{\alpha}))^{2}\int_{0}^{t}\|((t+\delta-s)^{\alpha-1}-(t-s)^{\alpha-1})h(s)\|^{2}_{L_{2}^{0}}ds\nonumber\\
&\quad+8Ht^{2H-1}\int_{0}^{t}\|(t-s)^{\alpha-1}(S_{\alpha}(t+\delta-s)-S_{\alpha}(t-s))h(s)\|^{2}_{L_{2}^{0}}ds.
\end{align}
For $p>\frac{1}{2\alpha-1}$, note that $\alpha\in (\frac{1}{2},1]$, we have
\begin{align}\label{h}
\int_{0}^{t}(t-s)^{2\alpha-2}\|h(s)\|^{2}_{L_{2}^{0}}ds&\leq \left(\int_{0}^{t}(t-s)^{\frac{(2\alpha-2)p}{p-1}}ds\right)^{\frac{p-1}{p}}\left(\int_{0}^{t}\|h(s)\|_{L_{2}^{0}}^{2p}ds\right)^{\frac{1}{p}}\nonumber\\
&\leq b^{\frac{(2\alpha-1)p-1}{p}}\left(\int_{0}^{b}\|h(s)\|_{L_{2}^{0}}^{2p}ds\right)^{\frac{1}{p}}\nonumber\\
&<\infty.
\end{align}
Similarly, we have
\begin{align}\label{s}
\int_{0}^{t}(t+\delta-s)^{2\alpha-2}\|h(s)\|^{2}_{L_{2}^{0}}ds<\infty.
\end{align}
Applying the dominated convergence theorem to (\ref{I1}), it follows that
\begin{align}\label{app}
I_{1}\rightarrow 0 \ \ \mbox{as} \ \ \delta\rightarrow 0.
\end{align}
For $I_{2}$, we have
\begin{align}\label{II}
I_{2}&\leq 8H  t^{2H-1}M^{2}(E_{\alpha,\alpha}(\omega(t+\delta)^{\alpha}))^{2}\int_{t}^{t+\delta}(t+\delta-s)^{2\alpha-2}\|h(s)\|_{L_{2}^{0}}^{2}ds\nonumber\\
&\leq 8H  t^{2H-1}M^{2}(E_{\alpha,\alpha}(\omega(t+\delta)^{\alpha}))^{2} \left(\int_{t}^{t+\delta}(t+\delta-s)^{\frac{(2\alpha-2)p}{p-1}}ds\right)^{\frac{p-1}{p}}\left(\int_{t}^{t+\delta}\|h(s)\|_{L_{2}^{0}}^{2p}ds\right)^{\frac{1}{p}}\nonumber\\
&\leq 8H  t^{2H-1}M^{2}(E_{\alpha,\alpha}(\omega(t+\delta)^{\alpha}))^{2}\delta^{\frac{(2\alpha-1)p-1}{p}}\left(\int_{0}^{b}\|h(s)\|_{L_{2}^{0}}^{2p}ds\right)^{\frac{1}{p}}.
\end{align}
By (\ref{II}), we have
\begin{align}\label{zero}
I_{2}\rightarrow 0 \ \ \mbox{as} \ \ \delta\rightarrow 0.
\end{align}
From (\ref{delta}), (\ref{app}), (\ref{zero}), it follows that
the stochastic integral term $\int_{0}^{t}(t-s)^{\alpha-1}S_{\alpha}(t-s)h(s)dB^{H}(s)$ belongs to the space $C([-r,b];L^{2}(\Omega;X))$.\\
Step 2: Set $x^{0}=0$ and construct a recurrence sequence process $\{x^{n}\}_{n\in \mathbb{N}}$ as
\begin{equation}\label{pro}
\left \{\begin{aligned}
&x^{n}(t)=T_{\alpha}(t)\phi(0)+\int_{0}^{t}(t-s)^{\alpha-1}S_{\alpha}(t-s)f(s,x^{n-1}_{s})ds\\
&\quad+\int_{0}^{t}(t-s)^{\alpha-1}S_{\alpha}(t-s)h(s)dB^{H}(s),\ t\in [0,b],\\
& x^{n}(t)=\phi(t),\ t\in[-r,0].
\end{aligned}\right.
\end{equation}
We will prove $x^{n}(t)\in C([-r,b];L^{2}(\Omega;X))$. Assume that $x^{n-1}(t)\in C([-r,b];L^{2}(\Omega;X))$. For sufficiently small $\delta>0$,
\begin{align}\label{I2}
\mathbb{E}\|x^{n}(t+\delta)-x^{n}(t)\|^{2}&\leq 2\mathbb{E}\|\int_{0}^{t}((t+\delta-s)^{\alpha-1}S_{\alpha}(t+\delta-s)-(t-s)^{\alpha-1}S_{\alpha}(t-s))f(s,x^{n-1}_{s})ds\|^{2}\nonumber\\
&\quad+2\mathbb{E}\|\int_{t}^{t+\delta}(t+\delta-s)^{\alpha-1}S_{\alpha}(t+\delta-s)f(s,x^{n-1}_{s})ds\|^{2}\nonumber\\
& =I_{3}+I_{4}.
\end{align}
For $I_{3}$, we have
\begin{align}\label{I3}
I_{3}&=2\mathbb{E}\|\int_{0}^{t}((t+\delta-s)^{\alpha-1}-(t-s)^{\alpha-1})S_{\alpha}(t+\delta-s)f(s,x^{n-1}_{s})ds\nonumber\\
&\quad+\int_{0}^{t}(t-s)^{\alpha-1}(S_{\alpha}(t+\delta-s)-S_{\alpha}(t-s))f(s,x^{n-1}_{s})ds\|^{2}\nonumber\\
&\leq 4\mathbb{E}\|\int_{0}^{t}((t+\delta-s)^{\alpha-1}-(t-s)^{\alpha-1})S_{\alpha}(t+\delta-s)f(s,x^{n-1}_{s})ds\|^{2}\nonumber\\
&\quad+4\mathbb{E}\|\int_{0}^{t}(t-s)^{\alpha-1}(S_{\alpha}(t+\delta-s)-S_{\alpha}(t-s))f(s,x^{n-1}_{s})ds\|^{2}\nonumber\\
&\leq 4M^{2}(E_{\alpha,\alpha}(\omega(t+\delta)^{\alpha}))^{2}\mathbb{E}\left(\int_{0}^{t}\|((t+\delta-s)^{\alpha-1}-(t-s)^{\alpha-1})f(s,x^{n-1}_{s})\|ds\right)^{2}\nonumber\\
&\quad+4\mathbb{E}\left(\int_{0}^{t}\|(t-s)^{\alpha-1}(S_{\alpha}(t+\delta-s)-S_{\alpha}(t-s))f(s,x^{n-1}_{s})\|ds\right)^{2}\nonumber\\
&\leq4M^{2}(E_{\alpha,\alpha}(\omega(t+\delta)^{\alpha}))^{2}\left(\int_{0}^{t}((t+\delta-s)^{\alpha-1}-(t-s)^{\alpha-1})^{2}ds\right)\mathbb{E}\int_{0}^{t}\|f(s,x^{n-1}_{s})\|^{2}ds\nonumber\\
&\quad+4\int_{0}^{t}(t-s)^{2\alpha-2}\|S_{\alpha}(t+\delta-s)-S_{\alpha}(t-s)\|^{2}ds\mathbb{E}\int_{0}^{t}\|f(s,x^{n-1}_{s})\|^{2}ds.
\end{align}
By $(f_{2}), \ (f_{3})$,
\begin{align}\label{expect}
\mathbb{E}\int_{0}^{t}\|f(s,x^{n-1}_{s})\|^{2}ds&=\mathbb{E}\int_{0}^{t}\|f(s,x^{n-1}_{s})-f(s,0)+f(s,0)\|^{2}ds\nonumber\\
&\leq 2\mathbb{E}\int_{0}^{t}\|f(s,x^{n-1}_{s})-f(s,0)\|^{2}ds+2\int_{0}^{t}\|f(s,0)\|^{2}ds\nonumber\\
&\leq 2C_{f}\mathbb{E}\int_{-r}^{t}\|x^{n-1}(s)\|^{2}ds+2\int_{0}^{t}\|f(s,0)\|^{2}ds\nonumber\\
&<\infty,
\end{align}
then applying the dominated convergence theorem to (\ref{I3}) to obtain
\begin{align}\label{dom}
I_{3}\rightarrow 0 \ \ \mbox{as} \ \ \delta\rightarrow 0.
\end{align}
For $I_{4}$, we get
\begin{align}\label{second}
I_{4}&\leq 2\mathbb{E}\left(\int_{t}^{t+\delta}\|(t+\delta-s)^{\alpha-1}S_{\alpha}(t+\delta-s)f(s,x^{n-1}_{s})\|ds\right)^{2}\nonumber\\
&\leq 2M^{2}(E_{\alpha,\alpha}(\omega(t+\delta)^{\alpha}))^{2}\int_{t}^{t+\delta}(t+\delta-s)^{2\alpha-2}ds\mathbb{E}\int_{0}^{t}\|f(s,x^{n-1}_{s})\|^{2}ds.
\end{align}
Hence
\begin{align}\label{sum}
I_{4}\rightarrow 0 \ \ \mbox{as} \ \ \delta\rightarrow 0.
\end{align}
By (\ref{I2}), (\ref{dom}), (\ref{sum}),
\begin{align*}
\mathbb{E}\|x^{n}(t+\delta)-x^{n}(t)\|^{2}\rightarrow 0 \ \ \mbox{as} \ \ \delta\rightarrow 0.
\end{align*}
Step 3: We will prove that $\{x^{n}\}$ is a Cauchy sequence. For $t\in[-r,b]$, since $x^{n}(t)=x^{n-1}(t)$  for $t\in[-r,0]$, we obtain
\begin{align}\label{u}
\mathbb{E}\|x^{n+1}(t)-x^{n}(t)\|^{2}&\leq\mathbb{E}\left(\int_{0}^{t}(t-s)^{\alpha-1}\|S_{\alpha}(t-s)(f(s,x_{s}^{n})-f(s,x_{s}^{n-1}))\|ds\right)^{2}\nonumber\\
&\leq M^{2}(E_{\alpha,\alpha}(\omega t^{\alpha}))^{2} \mathbb{E}\left(\int_{0}^{t}(t-s)^{\alpha-1}\|f(s,x_{s}^{n})-f(s,x_{s}^{n-1})\|ds\right)^{2}\nonumber\\
&\leq M^{2}(E_{\alpha,\alpha}(\omega t^{\alpha}))^{2}\int_{0}^{t}(t-s)^{2\alpha-2}ds\mathbb{E}\int_{0}^{t}\|f(s,x_{s}^{n})-f(s,x_{s}^{n-1})\|^{2}ds\nonumber\\
&\leq\frac{M^{2}(E_{\alpha,\alpha}(\omega t^{\alpha}))^{2}C_{f}t^{2\alpha-1}}{(2\alpha-1)}\int_{0}^{t}\mathbb{E}\|x^{n}(s)-x^{n-1}(s)\|^{2}ds.
\end{align}
Set $z^{n}(t)=\sup_{0\leq \theta\leq t}\mathbb{E}\|x^{n+1}(\theta)-x^{n}(\theta)\|^{2}$. By (\ref{u}), we have
\begin{align*}
z^{n}(t)\leq \frac{M^{2}(E_{\alpha,\alpha}(\omega t^{\alpha}))^{2}C_{f}t^{2\alpha-1}}{(2\alpha-1)}\int_{0}^{t}z^{n-1}(s)ds,\ n\in\mathbb{N},\ \ n\geq2.
\end{align*}
Proceeding by induction we obtain
\begin{align*}
z^{n}(t)\leq \frac{(\frac{M^{2}(E_{\alpha,\alpha}(\omega t^{\alpha}))^{2}C_{f}b^{2\alpha}}{2\alpha-1})^{n-1}}{(n-1)!}z^{1}(t),\  n\in\mathbb{N},\ \ n\geq2, \ t \in [0,b].
\end{align*}
This yields that $\{x^{n}\}_{n\in \mathbb{N}}$ is a Cauchy sequence in $C([-r,b];L^{2}(\Omega;X))$.\\
Step 4: We will show that the limit $y$ of $\{x^{n}\}_{n\in \mathbb{N}}$ is a solution of (\ref{fbm}). By $(f_{2})$ and (\ref{pro}),
\begin{align}\label{last}
&\mathbb{E}\|\int_{0}^{t}(t-s)^{\alpha-1}S_{\alpha}(t-s)f(s,x_{s}^{n-1})ds-\int_{0}^{t}(t-s)^{\alpha-1}S_{\alpha}(t-s)f(s,y_{s})ds\|^{2}\nonumber\\
&\leq \frac{M^{2}(E_{\alpha,\alpha}(\omega t^{\alpha}))^{2}C_{f}t^{2\alpha-1}}{(2\alpha-1)}\int_{0}^{t}\mathbb{E}\|x^{n-1}(s)-y(s)\|^{2}ds.
\end{align}
By (\ref{last}), (\ref{solution}), (\ref{pro}), we see that $y$ is the unique mild solution of (\ref{fbm}).  \ \ \ \ $\Box$\\
\textbf{Remark 3.4.} If we consider the special case that $\alpha=1$, then the problem (\ref{fbm}) reduces to the problem (1.1) in \cite{Caraballo}. Therefore, the existence result of mild solutions in this paper generalizes the corresponding existence result in \cite{Caraballo}.

\section{An Example}
Consider the following fractional partial differential equation
\begin{equation}\label{cauchy}
 \left\{\begin{aligned}
&d[J_{t}^{1-\alpha}((\xi(t,z)-\phi(0,z))]=\frac{\partial^{2}
\xi(t,z)}{\partial z^{2}}dt+\mu(t,\xi(t-r,z))dt+\gamma(t,z)dB^{H}(t), \  \ t\in[0,T], z\in [0,\pi],\\
&\xi(t,0)=\xi(t,\pi)=0,\ t\in(0,T),\\
&\xi(t,z)=\phi(t,z), t\in[-r,0],z\in[0,\pi],
\end{aligned}\right.
\end{equation}
where $\alpha\in (\frac{1}{2},1)$, $r>0$, $B^{H}(t)$ is a fractional Brownian motion with Hurst parameter $H\in(\frac{1}{2},1)$, $\phi:[-r,0]\times [0,\pi]\rightarrow R$ is a given function.

Let $X=L^{2}[0,\pi]$, $A=\frac{\partial^{2}}{\partial z^{2}}$,
$D(A)=\{\tau\in X, \tau, \tau' \  \mbox{are absolutely continuous},
\tau''\in X, \tau(0)=\tau(\pi)=0\}$.
Since $e_{n}=\sqrt{2/\pi}\sin nz, n\in N$, is the orthonormal system of eigenvectors of $A$, we see that $A\tau=-\Sigma_{n=1}^{\infty}n^{2}\langle\tau,e_{n}\rangle e_{n}$, $\tau\in D(A)$. It is well-known that $A$ is the generator of a strongly continuous semigroup of bounded linear operators $T(t)\tau=\sum_{n=1}^{\infty}e^{-n^{2}t}\langle\tau,e_{n}\rangle e_{n}$, $\tau\in X, \ t\geq 0$, and $\|T(t)\|\leq e^{-t}$. From Remark 3.2, it follows that $\|S_{\alpha}(t)\|\leq E_{\alpha,\alpha}(-t^{\alpha})$.

Set
$$\xi(t)z=\xi(t,z),\ t\in[0,b], \ z\in[0,\pi],$$
$$f(t,\phi)z=\mu(t,\phi(\theta,z)),\ \theta\in[-r,0], \ z\in[0,\pi],$$
$$\phi(\theta)z=\phi(\theta,z),\ \theta\in[-r,0], \ z\in[0,\pi],$$
$$h(t)z=\gamma(t,z),\ t\in[0,b], \ z\in[0,\pi],$$
then problem (\ref{fbm}) is the abstract version of problem (\ref{cauchy}). We can choose suitable functions $\mu,\ \gamma$ such that the conditions $(f_{1})$, $(f_{2})$, $(f_{3})$ and $(h_{1})$ are satisfied. Then by Theorem 3.3, the problem (\ref{cauchy}) has a unique mild solution.
\section{Conclusion}
In this paper, we obtain the existence and uniqueness of mild solutions to a class of stochastic fractional functional differential equations driven by fractional Brownian motion. In fact, it is interesting to investigate the longtime behavior and the regularity of  mild solutions in the future work. On the other hand, one can consider the qualitative behaviour of such equations, for example, attractors, invariant measures.

\end{document}